\documentclass[12pt,a4paper]{article}
\usepackage[latin1]{inputenc}
\usepackage{amsmath}
\usepackage{amsfonts}
\usepackage{amssymb}
\usepackage{fancyhdr}
\author{Jussi I. Tyhtila}
\title{Dispersion measure for symmetric, stable probability distributions}
\newtheorem{maaritelma}{Definition}
\newtheorem{teoreema}{Theorem}

\newtheorem{korollaari}{Corollary}

\begin{document}
\maketitle

\begin{abstract}
In this sudy, a dispersion measure for stable distributions is proposed by measuring the 'average curvature' of the characteristic function of a stable random variable. The 'average curvature' is to be called the \textit{characteristic curvature} and it resembles the concept of \textit{Fisher information}. The intellectual motivation comes from the \textit{uncertainty principle for Fourier transform pairs}. Characteristic curvature reduces to that of standard deviation in the Gaussian case and possesses great many analytical features, for example the scaling behaviour of the type $t^{\frac{1}{\alpha}}$ and its explicit connection to the Euler Gamma function. It is also directly related to the Hurst exponent $H$ and the fractal dimension $2-H$. 
\end{abstract}

\section{Introduction}
Stable distributions is an interesting and important class of probability distributions. They were discovered explicitly by Paul L\'{e}vy in 1925 \cite{lk}. They possess many interesting properties, most importantly they are by definition invariant under addition, up to a scale. Noteworthly they have power-law type of decay and therefore they are an excellent model for modelling many natural phenomena, such as earthquakes, financial returns, and a multitude of social phenomena such as size distributions of cities and firms \cite{scaling}. The major problem concerning them is that they have an infinite variance \cite{GK} and therefore their practical applicability is somewhat limited. Also they generally do not possess a density expressible in an analytic form. This study proposes a dispersion measure for them, drawing ideas from Fisher information, differential geometry and most importantly, the uncertainty principle for Fourier transform pairs \cite{Weyl}. The study begins with a brief discussion on characteristic functions and their relation to Fourier transforms and their properties, proceeds to a brief presentation of stable distributions and accumulates in defining a concept of \textit{characteristic curvature}, which is proposed as a suitable measure of dispersion for class of stable distributions.

\section{Characteristic functions of random variables and Fourier transforms of probability measures}
The following presentation will follow loosely \cite{GK}. In mathematics it is common to obtain solutions to given problems by transforming the problem in some other space, solve the problem there, and then transform the solution back to the original space. One of the most important transforms is the \textit{Fourier transform}. It is an integral transform that transforms some function, or measure in so called \textit{frequency space}. Fourier transform for a measure on $\mathbb{R}$ is defined as:
\begin{maaritelma}
Let $F(x)$ be a measure on measurable space $(\mathbb{R},\mathcal{B}_{\mathbb{R}})$, where $\mathcal{B}_{\mathbb{R}}$ is the Borel $\sigma$-algebra on $\mathbb{R}$, then it has a Fourier transform defined as:
\begin{equation}
\widehat{F(u)}=\int_{-\infty}^{\infty}e^{iux}dF(x)
\end{equation}
where $e^{iux}=\cos(ux)+i\sin(ux)$. Moreover, if there exists a function $f(x)$ such that dF(x)=f(x)dx, then it is equivalent to
\begin{equation}
\widehat{f(u)}=\int_{-\infty}^{\infty}e^{iux}f(x)dx
\end{equation}

\end{maaritelma}
\begin{maaritelma}
Given probability distribution $F(X\le x)$, $x\in \mathbb{R}$, the characteristic function of the random variable $X$ is defined as:
\begin{equation}
\varphi_X(u)=\mathbb{E}(e^{iuX})=\int_{-\infty}^{\infty}e^{iux}dF(X\le x)
\end{equation}
Since $\vert e^{iux}\vert =1$, $\varphi_X(u)$ is defined for all $u\in \mathbb{R}$ for every $F(X\le x)$.
\end{maaritelma}
\begin{maaritelma}
If a probability distribution $F(X\le x)$ has a density $f$, its characteristic function is:
\begin{equation}
\varphi_X(u)=\mathbb{E}(e^{iuX})=\int_{-\infty}^{\infty}e^{iux}f(x)dx
\end{equation}
That is, for probability distributions with density $f(x)$, the characteristic function is just the Fourier transform of the density $f(x)$.
\end{maaritelma}

We proceed with some properties of characteristic functions.
\begin{teoreema}
A characteristic function is uniformly continuous on $\mathbb{R}$ and satisfies the conditions
\begin{equation}
\varphi_X(0)=1 \; , \vert \varphi_X(u) \vert \le 1 \; (-\infty <u<\infty)
\end{equation}
Proof: see \cite{GK}
\end{teoreema}
\begin{teoreema}
Given random variables $X$ and $Y=aX+b$, where $a$ and $b$ are constants, then the characteristic functions of the random variables $X$ and $Y$ are connected by the equation
\begin{equation}
\varphi_Y(u)=\varphi_X(au)e^{ibu}
\end{equation}
and if $a>0$, then their distribution functions satisfy the relation
\begin{equation}
F(Y\le x)=F(X\le \frac{x-b}{a})
\end{equation}
\end{teoreema}
Proof. We note that
$$\varphi_Y(u)=\mathbb{E}(e^{iuY})=\mathbb{E}(e^{iu(aX+b)})=e^{iub}\mathbb{E}(e^{iuaX})=e^{iub}\varphi_X(au)$$
and for $a>0$
$$F(Y\le x)=F(aX+b\le x)=F(X\le \frac{x-b}{a})$$
QED.
\begin{teoreema}
The characteristic function of a random variable $X$ is real if and only if the distribution function of the random variable $X$ is symmetrical, that is for every x the following equation holds:
$$F(X\le x)+F(X\le -x)=1$$
\end{teoreema}
Proof. See \cite{GK}.
The following theorem is of great importance:
\begin{teoreema}
(Convolution theorem) The characteristic function of the sum of two independent random variables is the product of the characteristic functions of the summands. That is,
\begin{equation}
\varphi_{X+Y}(u)=\varphi_X(u)\varphi_Y(u) 
\end{equation}
\end{teoreema}
Proof. Given X and Y independent, then $e^{iuX}$ and $e^{iuY}$ are also independent.Therefore,
$$\varphi_{X+Y}(u)=\mathbb{E}(e^{iu(X+Y)})=\mathbb{E}(e^{iuX}e^{iuY})=\mathbb{E}(e^{iuX})\mathbb{E}(e^{iuY})=\varphi_X(u)\varphi_Y(u)$$
As expectation of the product of independent random variables is the product of the expectations, QED.

\section{Stable probability distributions} 
In this section, we introduce stable distributions and explore their most important properties. We shall follow loosely \cite{GK}. 
\begin{maaritelma}
The distribution function $F(X\le x)$ is called stable if to every $a_1>0$, $b_1$, $a_2>0$, $b_2$ there correspond constants $a>0$ and $b$ such that the equation 
\begin{equation}
F(X\le a_1x+b_1)\ast F(X\le a_2x+b_2)=F(X\le ax+b) 
\end{equation}
holds. 
\end{maaritelma}
Equivalently, we can define stable distributions as:
\begin{maaritelma}
The distribution function $F(X\le x)$ is called stable if to every $a_1>0$, $b_1$, $a_2>0$, $b_2$ there correspond constants $a>0$ and $b$ such that the equation 
\begin{equation}
F(\frac{X-b_1}{a_1}\le x)\ast F(\frac{X-b_2}{a_2}\le x)=F(\frac{X-b}{a}\le x) 
\end{equation}
holds. 
\end{maaritelma}
or, still:
\begin{maaritelma}
The distribution function $F(X\le x)$ is called stable if to every $a_1>0$, $b_1$, $a_2>0$, $b_2$ there correspond constants $a>0$ and $b$ such that the equation 
\begin{equation}
\varphi_{\frac{X-b}{a}}(u)=\varphi_{\frac{X-b_2}{a_2}}(u)\varphi_{\frac{X-b_1}{a_1}}(u)
\end{equation}
holds.
\end{maaritelma}

In words, a stable distribution $F$ is such, that given random variable $X$ and two affine transformations of it distributed with the same law $F$, then the sum of those transformed random variables is distributed according to some other affine transformation of $X$ with the same law $F$.

\subsection{Canonical representation of stable laws}
Next theorem gives an explicit function class for the characteristic functions of stable distributions.
\begin{teoreema}
In order that the distribution function $F(X\le x)$ is stable, it is necessary and sufficient that its characteristic function is:
\begin{equation}
\varphi_X(u)=e^{i\gamma u-c\vert u\vert^{\alpha}(1+i\beta\frac{u}{\vert u\vert}\omega(u,\alpha))}
\end{equation}
where $\alpha ,\beta ,\gamma ,c>0$ are constants($\gamma \in \mathbb{R}$, $-1\le \beta \le 1$, $0<\alpha \le 2$ and $\omega (u,\alpha)=tan(\frac{\pi \alpha}{2})$ if $\alpha \neq 1$ and $\omega (u,\alpha)=\frac{2}{\pi}\ln \vert u \vert$ if $\alpha =1$
\end{teoreema}
The proof is long and can be found in \cite{GK}. Following is important for our purpose of this study:
\begin{korollaari}
A stable distribution function $F(X\le x)$ that is symmetric around zero\footnote{In this study, the expressions \textit{symmetric around zero} and \textit{symmetric} are considered synonymous}, has the characteristic function
\begin{equation}
\varphi_X(u)=e^{-c\vert u \vert ^{\alpha}}
\end{equation}
\end{korollaari}
This results from the properties of characteristic functions discussed above.

\subsection{Sums of symmetric, stable random variables}
Now we examine the summing of independent, identically distributed symmetric, stable random variables. As the study of stochastic processes and therefore stochastic modelling is based on the properties of sums of random variables, in this chapter we study them in the case of symmetric, stable distributions. \\

Consider the set of independent, identically distributed random variables $X_i$ $i\in(1,2, ...,t)$ with $X_i$ having the characteristic function 
\begin{equation}
\varphi_{X_i}(u)=e^{-c\vert u\vert^{\alpha}}\; \forall i
\end{equation}
Let us consider the sum of such random variables:
\begin{equation}
S_t=\sum_{i=1}^{t}X_i=X_1+X_2+\cdots +X_t
\end{equation}
We wish to determine the characteristic function of $S_t$. As the convolution theorem states that the characteristic function of the sum of iid. random variables is the product of respective characteristic functions, we shall have:
\begin{equation}
\varphi_{S_t}(u)=\prod_{i=1}^{t}\varphi_{X_i}(u)
\end{equation}
Substituting,
\begin{equation}
\varphi_{S_t}(u)=(e^{-c\vert u\vert^{\alpha}})^t=e^{-tc\vert u\vert^{\alpha}}
\end{equation}
If we set $C=C(t)=ct$, we have:
\begin{equation}
\varphi_{S_t}(u)=(e^{-c\vert u\vert^{\alpha}})^t=e^{-C(t)\vert u\vert^{\alpha}}
\end{equation}
That is, the sum of $t$ iid. random variables $X_i$ has the same characteristic function as $X_i$ up to a linearly growing scaling factor $C$

\subsubsection{Example: The Gaussian case $N(0,\sigma ^2)$}
Consider random variable $X_i$ with characteristic function:
\begin{equation}
\varphi_{X_i}(u)=e^{-\frac{1}{2}\sigma ^2\vert u\vert^2}
\end{equation}
Now we have for the sum $S_t$ the characteristic function:
\begin{equation}
\varphi_{S_t}(u)=e^{-\frac{1}{2}\sigma ^2t\vert u\vert^{2}}
\end{equation}
If we set, 
\begin{equation}
\sigma^2(t)=\sigma^2t
\end{equation}
or equivalently:
\begin{equation}
\sigma(t)=\sigma \sqrt{t}
\end{equation}
we will have:
\begin{equation}
\varphi_{S_t}(u)=e^{-\frac{1}{2}\sigma ^2(t)\vert u\vert^{2}}
\end{equation}
As we can see, the variance of the sum grows linearly as a function of the number of summands. This kind of behaviour for sums of iid. random variables, that is the linear growing of variance as a function of time (summands) is called \textit{normal diffusion}. If one solves eg. the heat equation, one obtains a time-dependent normal distribution that has a linearly growing variance as a function of time \cite{pde}.

\subsection{Scaling of symmetric, stable distributions}
Next we concentrate on the scaling properties of stable distributions. As by definition the sum of $t$ stable iid. random variables has the same distribution up to \textit{scale}, the natural question is: \textit{what is that scale and how doest it evolve when the number of summands is increased}? We answer this question straightforwardly by presenting the following theorem: 
\begin{teoreema}
(Scaling theorem) Given symmetric, stable random variable $X$ with characteristic function $\varphi_X(u)=e^{-c\vert u\vert^{\alpha}}$ and distribution function $F(X\le x)$. Then the sum $S_t$ of $t$ iid. random variables $X$ has the distribution:
\begin{equation}
F(t^{\frac{1}{\alpha}}X\le x)
\end{equation}
or equivalently
\begin{equation}
F(X\le \frac{x}{t^{\frac{1}{\alpha}}})
\end{equation}
\end{teoreema}
Proof. We have the characteristic function:
\begin{equation}
\varphi_X(u)=e^{-c\vert u\vert^{\alpha}}
\end{equation}
Suppose we scale the random variable $X$ such that $X\longrightarrow sX$, where $s>0$. The characteristic function of the scaled random variable $sX$ will be:
\begin{equation}
\varphi_{sX}(u)=e^{-c\vert su\vert^{\alpha}}
\end{equation}
Which is equivalent to:
\begin{equation}
\varphi_{sX}(u)=e^{-cs^{\alpha}\vert u\vert^{\alpha}}
\end{equation}
In the previous section we showed, that the sum of $t$ independent, identically distributed symmetric, stable random variables $X$ has the characteristic function:
\begin{equation}
\varphi_{S_t}(u)=e^{-tc\vert u\vert^{\alpha}}
\end{equation}
by comparing, we note that these two expressions are the same, if $t=s^{\alpha}$ or $s=t^{\frac{1}{\alpha}}$. Hence if we scale the random variable $X$ by a factor $t^{\frac{1}{\alpha}}$, the scaled random variable $t^{\frac{1}{\alpha}}X$ has exactly the same distribution as the sum of $t$ iid. random variables $X$. In distribution function notation, the sum of $t$ independent, identically distributed symmetric, stable random variables $X$ has the distribution:
\begin{equation}
F(t^{\frac{1}{\alpha}}X\le x)
\end{equation}
or equivalently:
\begin{equation}
F(X\le \frac{x}{t^{\frac{1}{\alpha}}})
\end{equation}
where $F(X\le x)$ is the distribution of $X$. QED.
So the probability distribution spreads out proportionally to  $t^{\frac{1}{\alpha}}$ when number of summands $t$ is increased. This is the fractal property of sums of stable random variables, the sum distribution looks exactly the same as the original, it is just zoomed by a factor $t^{\frac{1}{\alpha}}$. This scaling property will come up later on in an interesting context.

\subsection{Some remarks on the scaling properties}
As we have seen, stable distributions obey a scaling law of the type $t^{\frac{1}{\alpha}}$ for sums of $t$ iid. symmetric, stable random variables. We have also seen that the scale parameter $c>0$ in the characteristic function grows linearly as a function of the number of summands (or time, if you will). These properties \textit{do} tell us a whole lot about the sum behaviour of stable random variables. In the Gaussian case one can directly state how the standard deviation grows over time if we consider the sum as a 1-dimensional random walk. This is extremely useful, as for example in the context of mathematical finance, standard deviation of the underlying price process is called \textit{volatility}. However, as the Gaussian case is the only case where variance is finite, we can not give predictions of the growth of standard deviation over time in other cases when $\alpha<2$. Basically we face the problem that we do know the time-evolution of the dispersion of the sum-distribution, but we do not know how to measure that dispersion in general. This is the topic of the next chapter.

\section{Characteristic curvature}

\subsection{Motivation}
We next proceed to propose a measure for dispersion for symmetric around zero stable distributions with $\alpha \in [1,2]$.
The idea is to think of a measure that somehow characterises the dispersive properties of the characteristic functions. We know that characteristic functions determine probability distributions completely \cite{GK}, whether the moments are finite or not. So the general dispersive properties of different probability distributions must be somehow incorporated in the properties of characteristic functions.We know from Fourier analysis, that the more spread the original function, the more concentrated its Fourier transform and vice versa \cite{Weyl}. Because characteristic function is the Fourier transform of the probability measure, we can apply this idea here. It is therefore natural to think that if the probability distribution is greatly dispersed in space, the characteristic function should be condensed around the origin in the frequency space. This concentration of characteristic functions around the origin is intuitively analogous to that the second derivative in its absolute value should be relatively large around the origin, because the characteristic function satisfies $\varphi_X(0)=1$ and it should decay fast when moved from the origin into either direction. Hower, because the characteristic function here is symmetric around the origin, we can constrain to analyse the positive part of $u$. Therefore the negative of the second derivative should be relatively large for small $u$. To measure this concentration analytically, let us consider the weighted average of the negative second derivatives of $\ln$-transform of $\varphi_X(u)$ (assuming $\varphi_X(u)$ monotonically decreasing $\forall u\ge 0$) in the continuum limit, weights given by the characteristic function. That is, the integral transform that we shall call as
the \textit{characteristic curvature}:
\begin{maaritelma}
The characteristic curvature for symmetric around zero probability distribution $F(X\le x)$ is
\begin{equation}
\Phi=\frac{2}{\sqrt{2\pi}}\int_{0}^{\infty}(-\frac{\partial^2 \ln \varphi_X(u)}{\partial u^2})\varphi_X(u)du
\end{equation}
where $\varphi_X(u)$ is the characteristic function of the distribution and $u\ge 0$.
\end{maaritelma}
The constant $\frac{2}{\sqrt{2\pi}}$ is suitable for fixing the proper scale, as we shall see later on.
Note that this integral resembles the definition of \textit{Fisher information}, differences being that now the likelihood function is the characteristic function, and the parameter is $u$. By performing the derivation and simplifying, the characteristic curvature can be formulated equivalently as:
\begin{equation}
\Phi=\frac{2}{\sqrt{2\pi}}\int_{0}^{\infty}( \frac{\varphi'_X(u)^2}{\varphi_X(u)}-\varphi''_X(u))du
\end{equation}
which resembles some curvature concepts in differential geometry. For the sake of comparison, the standard deviation of a probability distribution with characteristic function $\varphi_X(u)$ is given as:
\begin{equation}
\sigma=\sqrt{-(\varphi'_X(0))^2-\varphi''_X(0)}
\end{equation}
see eg. \cite{GK}.

\subsection{Characteristic curvature for symmetric, stable distributions}
\begin{teoreema}
Given symmetric around zero stable distribution with characteristic function $\varphi_{X}(u)=e^{-c\vert u \vert^\alpha}$ its characteristic curvature satisfies
\begin{equation}
\Phi(c;\alpha)=\frac{2}{\sqrt{2\pi}}\int_{0}^{\infty}\lbrace -\frac{\partial ^2 \ln \varphi_X(u)}{\partial u^2} \rbrace e^{-cu^\alpha}du
=\frac{2}{\sqrt{2\pi}}\alpha \Gamma(2-\frac{1}{\alpha})c^{\frac{1}{\alpha}}=
\frac{2}{H\sqrt{2\pi}}\Gamma(2-H)c^{H}\end{equation}
where $\Gamma(2-\frac{1}{\alpha})$ is the Euler Gamma function and $H$ is the Hurst exponent.
\end{teoreema}
Proof. \\
Let us begin with the characteristic function of a symmetric stable distribution for random variable X.
\begin{equation}
\varphi_{X}(u)=e^{-c\vert u \vert^\alpha} \; \; u\in \mathbb{R}, \; c>0, \; 0<\alpha\le 2
\end{equation}
Because the characteristic function is symmetric around zero, we consider the positive half-space representation of it:
\begin{equation}
\varphi_{X}(u)=e^{-cu ^\alpha} \; \; u\ge 0, \; c>0, \; 0<\alpha\le 2
\end{equation}
Taking logarithms, we will get
\begin{equation}
\ln \varphi_{X}(u)=-cu^\alpha
\end{equation}
Derivating twice with respect to $u$:
\begin{equation}
\frac{\partial^2 \ln \varphi_{X}(u) }{\partial u^2}=-\alpha (\alpha -1)cu^{\alpha -2}
\end{equation} 
Which is the same as
\begin{equation}
-\frac{\partial^2 \ln \varphi_{X}(u) }{\partial u^2}=\alpha (\alpha -1)cu^{\alpha -2}
\end{equation}
Substituting this into characteristic curvature, we shall have:
\begin{equation}
\Phi(c;\alpha)=\frac{2}{\sqrt{2\pi}}\int_{0}^{\infty}\lbrace \alpha (\alpha -1)cu^{\alpha -2}  \rbrace e^{-cu^\alpha}du
\end{equation}
Which is equivalent to:
\begin{equation}
\Phi(c;\alpha)=\frac{2}{\sqrt{2\pi}} \alpha (\alpha -1)c\int_{0}^{\infty}u^{\alpha -2}e^{-cu^\alpha}du
\end{equation}
Next we calculate the integral:
\begin{equation}
\int_{0}^{\infty}u^{\alpha -2}e^{-cu^\alpha}du
\end{equation}
Substituting $t=cu^\alpha \; \Leftrightarrow dt=c\alpha u^{\alpha -1}du \; \Leftrightarrow du=\frac{dt}{c\alpha u^{\alpha -1}}$ The integral will become:
\begin{equation}
\int_{0}^{\infty}\frac{u^{\alpha -2}e^{-t}}{c\alpha u^{\alpha -1}}dt
\end{equation}
which is
\begin{equation}
\frac{1}{c\alpha}\int_{0}^{\infty}u^{-1}e^{-t}dt
\end{equation}
Noting that $t=cu^{\alpha} \; \Leftrightarrow u=\frac{t^{\frac{1}{\alpha}}}{c^{\frac{1}{\alpha}}}$ the integral becomes
\begin{equation}
\frac{1}{c\alpha} \int_{0}^{\infty} \frac{ t^{-\frac{1}{\alpha} } }   {c^{-\frac{1}{\alpha}}}    e^{-t}dt
\end{equation}
Taking the constant outside the integral gives:
\begin{equation}
\frac{1}{\alpha cc^{-\frac{1}{\alpha}}} \int_{0}^{\infty}t^{-\frac{1}{\alpha}}e^{-t}dt
\end{equation}
By simplifying and adding zero to the exponent of $t$, we shall have:
\begin{equation}
\frac{1}{\alpha c^{1- \frac{1}{\alpha}}} \int_{0}^{\infty}t^{(1-\frac{1}{\alpha})-1}e^{-t}dt
\end{equation}
Recalling that the Euler Gamma function is defined as
\begin{equation}
\Gamma (z)=\int_{0}^{\infty}t^{z-1}e^{-t}dt \; z\in \mathbb{C}
\end{equation}
We will obtain the elegant result:
\begin{equation}
\int_{0}^{\infty}u^{\alpha -2}e^{-cu^\alpha}du=\frac{ c^ {\frac{1}{\alpha}-1}} {\alpha} \Gamma(1-\frac{1}{\alpha})
\end{equation}
Substituting this to the characteristic curvature, we will obtain:
\begin{equation}
\Phi(c;\alpha)=\frac{2}{\sqrt{2\pi}} \alpha (\alpha -1)c\frac{ c^ {\frac{1}{\alpha}-1}} {\alpha} \Gamma(1-\frac{1}{\alpha})
\end{equation}
Simplifying, it will become:
\begin{equation}
\Phi(c;\alpha)=\frac{2}{\sqrt{2\pi}}(\alpha -1)c^{\frac{1}{\alpha}} \Gamma(1-\frac{1}{\alpha})
\end{equation}
Because the Gamma function satisfies
\begin{equation}
\Gamma(z+1)=z\Gamma(z)
\end{equation}
We have
\begin{equation}
\Gamma(2-\frac{1}{\alpha})=(1-\frac{1}{\alpha})\Gamma(1-\frac{1}{\alpha}) \; \Leftrightarrow \Gamma(1-\frac{1}{\alpha})=\frac{\alpha}{\alpha -1}\Gamma(2-\frac{1}{\alpha})
\end{equation}
Substituting this in the characteristic curvature:
\begin{equation}
\Phi(c;\alpha)=\frac{2}{\sqrt{2\pi}}(\alpha -1)c^{\frac{1}{\alpha}}\frac{\alpha}{\alpha -1}\Gamma(2-\frac{1}{\alpha})=\frac{2}{\sqrt{2\pi}}\alpha c^{\frac{1}{\alpha}}\Gamma(2-\frac{1}{\alpha})
\end{equation}
We will finally get the most compact form for characteristic curvature:
\begin{equation}
\Phi(c;\alpha)=\frac{2}{\sqrt{2\pi}}\alpha \Gamma(2-\frac{1}{\alpha})c^{\frac{1}{\alpha}}
\end{equation}
QED

\subsection{Properties of characteristic curvature for stable distributions and some restrictions}
We have shown that the characteristic curvature for stable distributions satisfies:
$$
\Phi(c;\alpha)=\frac{2}{\sqrt{2\pi}}\alpha \Gamma(2-\frac{1}{\alpha})c^{\frac{1}{\alpha}}
$$
As the gamma function has singularities at $x=0,-1,-2,-3...$, characteristic curvature diverges when $\alpha$ approaches $\frac{1}{2}$. If we restrict that 
$1\le \alpha \le 2$, we will have following desirable properties for characteristic curvature:
\begin{teoreema}
Given $1\le \alpha \le 2$ and $\forall c>0$, the characteristic curvature $\Phi(\alpha ,c)$ satisfies the following properties: 
$$\begin{array}{cc}
(1)& \Phi(c;\alpha)>0 \\
(2)& \frac{\partial \Phi(c;\alpha) }{\partial c} > 0
\end{array}
$$
\end{teoreema}
Proof. (1) is trivial. We prove hence only (2). By derivating:
$$\frac{\partial \Phi}{\partial c}=\frac{2}{\sqrt{2\pi}}\Gamma(2-\frac{1}{\alpha})c^{\frac{1-\alpha}{\alpha}}> 0$$
QED.
In other words, characteristic curvature is strictly increasing in $c$, which is the scale factor, and is always positive.

\subsection{Characteristic curvature for Gaussian distribution}
Let us calculate the characteristic curvature for centered Normal distribution with variance $\sigma^2$.
The characteristic function is
$$\varphi_X(u)=e^{-\frac{1}{2}\sigma^2u^2}$$
Equivalently $\alpha=2$ and $c=\frac{1}{2}\sigma^2$. By substituting these values, we will get:
$$\Phi(2 ,\frac{1}{2}\sigma^2)=\frac{4}{\sqrt{2\pi}}\sqrt{\frac{1}{2}\sigma^2}\Gamma(\frac{3}{2})$$
Because $\Gamma(\frac{3}{2})=\frac{\sqrt{\pi}}{2}$, we will obtain 
$$\Phi(2 ,\frac{1}{2}\sigma^2)=\frac{4}{\sqrt{2\pi}}\sqrt{\frac{1}{2}\sigma^2}\frac{\sqrt{\pi}}{2} =\frac{2\sqrt{\pi}}{\sqrt{2\pi}}\sqrt{\frac{1}{2}\sigma^2}$$
Which simplifies to:
$$\Phi(2 ,\frac{1}{2}\sigma^2)=\sqrt{2}\sqrt{\frac{1}{2}}\sigma=\sigma$$
This is very fascinating result! The characteristic curvature of a centered Gaussian with variance $\sigma^2$ is just the standard deviation $\sigma$.
Note that also holds:
\begin{equation}
\Phi^2(2 ,\frac{1}{2}\sigma^2)=\sigma^2
\end{equation}
,which is obvious, of course.
\subsection{Characteristic curvature for Cauchy distribution}
Let us calculate the characteristic curvature for Cauchy distribution, that is $\alpha=1$ and $c$ is unspecified.
The characteristic function is:
$$\varphi_X(u)=e^{-c\vert u \vert}$$
Substituting these values we will have:
$$\Phi(1 ,c)=\frac{2c}{\sqrt{2\pi}}\Gamma(1)$$
Because $\Gamma(1)=1$, we will have simply:
$$\Phi(1 ,c)=\frac{2c}{\sqrt{2\pi}}$$

\subsection{Sums of symmetric, stable random variables and characteristic curvature }
We proceed to consider the characteristic curvature for sums of independent, identically distributed symmetric, stable random variables $X_i$ with characteristic function $\varphi_{X_i}(u)=e^{-\vert u\vert^{\alpha}}$ $\forall i\in (1,2,3,...,t)$. We recall that the sum
\begin{equation}
S_t=\sum_{i=1}^{t}X_i=X_1+X_2+\cdots +X_t
\end{equation}
has the characteristic function
\begin{equation}
\varphi_{S_t}(u)=e^{-tc\vert u\vert^{\alpha}}
\end{equation}
We calculate the respective characteristic curvature:
\begin{equation}
\Phi(tc;\alpha)=\frac{2}{\sqrt{2\pi}}\alpha \Gamma(2-\frac{1}{\alpha})(tc)^{\frac{1}{\alpha}}
\end{equation}
which becomes
\begin{equation}
\Phi(tc;\alpha)=\frac{2}{\sqrt{2\pi}}\alpha \Gamma(2-\frac{1}{\alpha})c^{\frac{1}{\alpha}}t^{\frac{1}{\alpha}}
\end{equation}
which equals
\begin{equation}
\Phi(tc;\alpha)=\Phi(c;\alpha)t^{\frac{1}{\alpha}}
\end{equation}
This is very elegant result! The characteristic curvature scales like $\sim t^{\frac{1}{\alpha}}$, which is exactly the same scaling behaviour that stable random variables possess in general. 

\subsubsection{Example: $X_i\sim N(0,\sigma^2)$}
Given that $X_i\sim N(0,\sigma^2)$ $\forall i$. That is, we have $c=\frac{1}{2}\sigma^2$, $\alpha =2$, so for the sum of $t$ such random variables, we have:
\begin{equation}
\Phi(t\frac{1}{2}\sigma^2;2)=\Phi(\frac{1}{2}\sigma^2;2)t^{\frac{1}{2}}
\end{equation}
We recall from earlier results that $\Phi(\frac{1}{2}\sigma^2;2)=\sigma$, therefore we obtain:
\begin{equation}
\Phi(t\frac{1}{2}\sigma^2;2)=\sigma \sqrt{t}
\end{equation}
Which is exactly the same result we had earlier for normal diffusive behaviour.
Note that the square of the characteristic curvature is additive in the Gaussian case , that is when it is the variance, as:
\begin{equation}
\Phi^2(tc;\alpha)=\Phi^2(c;\alpha)t^{\frac{2}{\alpha}}=\sigma^2t
\end{equation}

\subsubsection{Example: $X_i\sim Cauchy\;(c)$}
Let $X_i$ have a Cauchy distribution with some scale factor $c$. We have $\alpha =1$, so for the sum of $t$ such random variables, we have:
\begin{equation}
\Phi(tc;1)=\Phi(c;1)t^{\frac{1}{1}}
\end{equation}
recalling earlier results, we then have:
\begin{equation}
\Phi(tc;1)=\frac{2c}{\sqrt{2\pi}}t
\end{equation}
which is the defining characteristics of so called \textit{ballistic diffusion}.

\section{Conclusions}
We have proposed a dispersion measure for stable distributions with $1\le \alpha \le 2$. Characteristic curvature has very intuitively appealing roots and interesting connections to various branches of mathematics. The scaling property; the familiar results of superdiffusive characteristics and the reduction to standard deviation in the Gaussian case are interesting implications of the definition. The abundance of Gamma function might lead to other interesting connections of stable distributions and e.g. Riemann zeta function as they are intimately related to each other. Obviously, more research is vitally needed to establish general results for the applicability of characteristic curvature for general probability distributions. Also, one interesting field of applied research would be that of option pricing for stable price processes.


\begin{thebibliography}{99}
\bibitem{lk} A. Ya. Khintcine, Paul Levy \emph{Sur les lois stables, C. R. Acad. Sci. Paris 202, 374-376 (1936)}
\bibitem{pde} L. C. Evans \emph{Partial differential equations, American mathematical society, ISBN 0-8218-0772-2 (1998)}
\bibitem{GK} B. V. Gnedenko, A. N. Kolmogorov \emph{Limit distributions for sums of independent random variables, Addison-Wesley (1954)}
\bibitem{JP} J. Jacod, P. Protter \emph{Probability essentials, Springer ISBN 3-540-43871-8 (2004)}
\bibitem{Weyl} H. Weyl \emph{Theory of groups and quantum mechanics, Dover, NY (1950)}
\bibitem{scaling} B. B. Mandelbrot \emph{Fractals and scaling in finance, Springer ISBN 0-387-98363-5 (1997)}
\end{thebibliography}
\end{document}